\documentclass[11pt]{article}
\usepackage{amssymb}
\usepackage{amsmath}
\usepackage{amsthm}
\usepackage{amsfonts}
\usepackage{url}
\usepackage{hyperref}
\usepackage{breakurl}

\newcommand{\Otilde}{\makebox{$\widetilde O$}}

\def\GF{\mathrm{GF}}
\def\gcd{\mathrm{gcd}}

\newtheorem{thm}{Theorem}
\newtheorem{hypothesis}{Assumption}

\begin{document}
\bibliographystyle{plain}

\title{{The Great Trinomial Hunt}%
\thanks{\hbox{To appear in {\em Notices of the AMS.}}
\hbox{Copyright \copyright\ R. Brent and P. Zimmermann, 2009.}
\hspace*{\fill}rpb235}
}
\author{Richard P. Brent and Paul Zimmermann}
\date{20 October 2009}

\maketitle

\section{Introduction}

\thispagestyle{empty}
A {\em trinomial} is a polynomial in one variable with three nonzero terms, 
for example $P = 6x^7 + 3x^3 - 5$. If the coefficients of a polynomial $P$ 
(in this case $6, 3, -5$)
are in some ring or field $F$, we say that $P$ is a polynomial over $F$,
and write $P \in F[x]$.
The operations of addition and multiplication of polynomials in $F[x]$ are
defined in the usual way, with the operations on coefficients performed in
$F$.

Classically the most common cases are $F = {\mathbf Z}, {\mathbf Q},
{\mathbf R}$ or ${\mathbf C}$, respectively the integers, rationals, reals
or complex numbers.  However, polynomials over finite fields are also
important in applications.  We restrict our attention to
polynomials over the simplest finite field: the field $\GF(2)$ of two
elements, usually written as $0$ and~$1$.
The field operations
of addition and multiplication are defined as for integers modulo~$2$, so
$0+1=1$, $1+1=0$, $0 \times 1 = 0$, $1 \times 1 = 1$, etc.

An important consequence of the definitions is that, for polynomials
$P, Q \in GF(2)[x]$, we have 
\[(P+Q)^2 = P^2+Q^2\] 
because the ``cross term''
$2PQ$ vanishes. High school algebra would have been much easier if we had
used polynomials over $\GF(2)$ instead of over ${\mathbf R}$!

Trinomials over $\GF(2)$ are important in cryptography and random number
generation. To illustrate why this might be true, consider 
a sequence $(z_0, z_1, z_2, \ldots)$ satisfying the recurrence
\begin{equation}
z_n = z_{n-s} + z_{n-r} \bmod 2,	\label{eq:recurrence}
\end{equation}
where $r$ and $s$ are given positive integers, $r > s > 0$, and the
initial values $z_0, z_1, \ldots, z_{r-1}$ are also given. The recurrence
then defines all the remaining terms $z_r, z_{r+1}, \ldots$ in the sequence.

It is easy to build hardware to implement the
recurrence~(\ref{eq:recurrence}). All we need is a shift register
capable of storing $r$ bits, and a circuit capable of computing the addition
mod~$2$ (equivalently, the ``exclusive or'') of two bits separated by 
$r-s$ positions in the shift register and feeding the output back into
the shift register. 
%

The recurrence (\ref{eq:recurrence}) 
looks similar to the well-known Fibonacci recurrence
\[F_n = F_{n-1} + F_{n-2};\]
indeed the Fibonacci numbers mod~$2$ satisfy our recurrence with $r = 2$,
$s = 1$.  This gives a sequence $(0, 1, 1, 0, 1, 1, \ldots)$ with
period~$3$: not very interesting.  However, if we take $r$ larger we can get
much longer periods. 

The period can be as large as $2^r-1$, which
makes such sequences interesting as components in pseudo-random number
generators or stream ciphers. In fact, the period is $2^r-1$ if the 
initial values are not all zero and the associated trinomial
\[x^r + x^s + 1,\]
regarded as a polynomial over $\GF(2)$, is {\em primitive}. 
A primitive polynomial is one that is irreducible (it has 
no nontrivial factors), and satisfies an additional condition
given in the ``Mathematical Foundations'' section below.

A {\em Mersenne prime} is a prime of the form $2^r-1$. 
Such primes are named after
Marin Mersenne (1588--1648), who corresponded with many of the scholars of
his day, and in 1644 gave a list (not quite correct) of the Mersenne 
primes with $r \le 257$.

A {\em Mersenne exponent} is the exponent $r$ of a 
Mersenne prime $2^r-1$. A Mersenne exponent is necessarily prime, but not
conversely. For example, $11$ is not a Mersenne exponent because
$2^{11}-1 = 23\cdot 89$ is not prime.

The topic of this article is a search for primitive trinomials of large
degree~$r$, and its interplay with a search for large Mersenne primes.
First, we need to explain the connection between these two topics,
and briefly describe the GIMPS project. Then we describe the algorithms
used in our search, which can be split into two distinct periods,
``classical'' and ``modern''. Finally, we describe the results obtained
in the modern period.

\section{Mathematical Foundations}

As stated above, we consider polynomials over the finite field $\GF(2)$.
An \emph{irreducible polynomial} is a polynomial that is not divisible by
any non-trivial polynomial other than itself.
For example $x^5+x^2+1$ is
irreducible, but $x^5 + x + 1$ is not, since 
$x^5 + x + 1 = (x^2 + x + 1)(x^3 + x^2 + 1)$ in $\GF(2)[x]$.
We do not consider binomials $x^r + 1$, because they are divisible by
$x + 1$, and thus reducible for $r > 1$. 

An irreducible polynomial $P$ of
degree $r>1$ yields a representation of the finite field $\GF(2^r)$ of $2^r$
elements: any polynomial of degree less than $r$ represents an element, the
addition is polynomial addition, whose result still has degree less than
$r$, and the multiplication is defined modulo $P$: one first multiplies both 
inputs, and then reduces their product modulo $P$.
Thus $\GF(2^r) \simeq \GF(2)[x]/P(x)$.

An irreducible polynomial $P$ of degree $r > 0$ over $\GF(2)$ 
is said to be \emph{primitive} iff $P(x) \ne x$ and the residue classes
$x^k \bmod P,\, 0 \leq k < 2^r-1$, are distinct.
In order to check primitivity of an irreducible
polynomial P, it is only necessary to check that $x^k \ne 1 \bmod P$ for those 
$k$ that are maximal non-trivial divisors of $2^r-1$.
For example, $x^5 + x^2 + 1$ is primitive; $x^6 + x^3 + 1$ is
irreducible but not primitive, since $x^9 = 1 \bmod (x^6+x^3+1)$.
Here $9$ divides $2^6-1=63$ and is a maximal divisor as $63/9 = 7$ is prime.

We are interested in primitive polynomials because
$x$ is a generator of the multiplicative
group of the finite field $\GF(2)[x]/P(x)$ if $P(x)$ is primitive.

If $r$ is large and $2^r-1$ is not
prime, it can be difficult to test primitivity of a polynomial of degree
$r$, because we need to know the prime factors of $2^r-1$. Thanks to the
Cunningham project~\cite{Cunningham}, 
these are known for all $r < 887$, %
but not in general for larger $r$.
On the other hand, if $2^r-1$ is prime, 
then all irreducible polynomials of degree $r$ are
primitive. This is the reason why we consider degrees
$r$ that are Mersenne exponents.

\section{Starting the Search}

In the year 2000 the authors were communicating by email with each other and
with Samuli Larvala when the topic of efficient algorithms for testing
irreducibility or primitivity of trinomials over $\GF(2)$ arose. The first
author had been interested in this topic for many years because of the
application to pseudo-random number generators. %
Publication of
a paper by Kumada {\em et al.}~\cite{Kumada00}, describing a search for
primitive trinomials of degree $859\,433$ 
(a Mersenne exponent), prompted the
three of us to embark on a search for primitive trinomials of degree
$r$, for $r$ ranging over all known Mersenne exponents. 
At that time, the largest known
Mersenne exponents were $3\,021\,377$ and $6\,972\,593$. 
The existing programs took time proportional to $r^3$. Since
$(6972593/859433)^3 \approx 534$, 
and the computation by Kumada {\em et al.}
had taken three months on $19$ processors, 
it was quite a challenge.

\section{The GIMPS project}

GIMPS stands for \emph{Great Internet Mersenne Prime Search}. It is a
distributed computing project started by George Woltman, with home page
\url{www.mersenne.org}. The goal of GIMPS is to find new Mersenne
primes. %
As of December 2009, GIMPS has found $13$ new Mersenne primes in 13 years,
and has held the record of the largest known prime since
the discovery of $M_{35}$ in 1996.
Mersenne primes are usually numbered in increasing order of size: 
$M_1 = 2^2-1 = 3$, $M_2 = 2^3-1 = 7$,
$M_3 = 2^5-1 = 31$, $M_4 = 2^7-1 = 127$, \ldots, 
$M_{38} = 2^{6972593}-1$, etc.

Since GIMPS does not always find Mersenne primes in order,
there can be some uncertainty in numbering the largest known Mersenne
primes. We write $M'_n$ for the $n$-th Mersenne prime in order of 
discovery.
There are gaps in the search above $M_{39} = 2^{13466917}-1$.
Thus we can have $M'_n >  M'_{n+1}$ for $n > 39$. For example,
$M'_{45} = 2^{43112609} - 1$ was found before $M'_{46} = 2^{37156667} - 1$
and $M'_{47} = 2^{42643801}-1$.
At the time of writing this article, $47$ Mersenne
primes are known, and the largest is $M'_{45} = 2^{43112609} - 1$.

It is convenient to write $r_n$ for the exponent of $M_n$,
and $r'_n$ for the exponent of $M'_n$.  For example,
$r'_{45} = 43\,112\,609$.

\section{Swan's Theorem}

We state a useful theorem, known as Swan's theorem,
although the result was found much earlier by
Pellet~\cite{Pellet78} and Stickelberger~\cite{Stickelberger97}.
In fact, there are several theorems in Swan's paper~\cite{Swan62}.
We state a simplified version of Swan's Corollary~5.

\begin{thm} 				\label{thm:Swan}
Let $r > s > 0$, and assume $r+s$ is odd.
Then $T_{r,s}(x) = x^r + x^s + 1$ has an even number of irreducible
factors over $\GF(2)$ in the following cases:\\
 a) $r$ even, $r \ne 2s$, $rs/2 = 0$ or $1$ mod $4$.\\
 b) $r$ odd, $s$ not a divisor of $2r$, $r = \pm 3$ mod $8$.\\
 c) $r$ odd, $s$ a divisor of $2r$, $r = \pm 1$ mod $8$.\\
In all other cases $x^r + x^s + 1$ %
has an odd number of irreducible factors. 
\end{thm}
If both $r$ and $s$ are even, then $T_{r,s}$ is a square and has an even
number of irreducible factors.
If both $r$ and $s$ are odd, we can apply the theorem to 
the ``reciprocal polynomial''
$T_{r,r-s}(x) = x^rT(1/x) = x^r + x^{r-s} + 1$,
since $T_{r,s}(x)$ and $T_{r,r-s}(x)$ have the same number of
irreducible factors.

For $r$ an odd prime, and excluding the easily-checked
cases $s = 2$ or $r-2$,
case~(b) says that the trinomial has an
{even} number of irreducible factors, and hence must be {reducible},
if $r = \pm 3 \bmod 8$. Thus, we only need to consider those Mersenne
exponents with $r = \pm 1 \bmod 8$. Of the $14$ known Mersenne exponents
$r > 10^6$, only $8$ satisfy this condition. 

\section{Cost of the Basic Operations}

The basic operations that we need are squarings modulo the trinomial
$T = x^r + x^s + 1$, multiplications modulo $T$, and greatest common
divisors (GCDs) between
$T$ and a polynomial of degree less than $r$. We measure the cost
of these operations in terms of the number of bit or word-operations
required to implement them.
In $\GF(2)[x]$, squarings cost $O(r)$, due to the fact
that the square of $x^i + x^j$ is $x^{2i} + x^{2j}$.
The reduction modulo $T$ of a polynomial of degree less than $2r$ 
costs $O(r)$, due to the sparsity of $T$; thus modular squarings cost
$O(r)$.

Modular multiplications cost $O(M(r))$, where $M(r)$ is the cost of
multiplication of two polynomials of degree less than $r$ over $\GF(2)$;
the reduction modulo $T$ costs $O(r)$, so the multiplication
cost dominates the reduction cost.
The ``classical'' polynomial multiplication algorithm has
$M(r) = O(r^2)$, but an algorithm\footnote{This algorithm differs from the
Sch\"onhage-Strassen integer-multiplication algorithm, which does not work 
over $\GF(2)$. For details see~\cite{rpb232,Schonhage77}.} 
due to Sch\"onhage has $M(r) = O(r \log r \log \log r)$~\cite{Schonhage77}.

A GCD computation for polynomials of
degree bounded by $r$ costs $O(M(r) \log r)$ 
using a ``divide and conquer'' approach combined with Sch\"onhage's fast
polynomial multiplication.
The costs are summarized in Table~1.\\

\begin{table}[ht]	
\begin{center}
\begin{tabular}{|c|c|} \hline
modular squaring & $O(r)$ \\
modular product  & $O(M(r))$ \\
GCD              & $O(M(r) \log r)$ \\
\hline
\end{tabular}
\end{center}
\caption{Cost of the basic operations.}
\label{tab:table1}
\end{table}

\section{Testing Irreducibility}

Let ${\mathbf P}_r(x) = x^{2^r} -x$.  As was known to Gauss,
${\mathbf P}_r(x)$ is the product of all irreducible polynomials of
degree $d$, where $d$ runs over the divisors of $r$.  For example,
\[{\mathbf P}_3(x) = x(x+1)(x^3 + x + 1)(x^3 + x^2 + 1)\]
in $\GF(2)[x]$.  Here $x$ and $x+1$ are the irreducible polynomials of
degree~$1$, and the other factors are the irreducible polynomials of
degree~$3$.  Note that we can always write ``$+$'' instead of ``$-$'' when
working over $\GF(2)$, since $1 = -1$ (or, equivalently, $1+1 = 0$).

In particular, if $r$ is an odd prime, 
then a polynomial $P(x) \in \GF(2)[x]$ with
degree~$r$ is irreducible iff 
\begin{equation}
x^{2^r} = x \bmod P(x)\;.		\label{eq:irred-cond}
\end{equation}
(If $r$ is not prime, then~(\ref{eq:irred-cond}) is necessary but not
sufficient: we have to check a further condition to guarantee irreducibility,
see~\cite{vzGaGe99}.)

When $r$ is prime, equation~(\ref{eq:irred-cond})
gives a simple test for irreducibility (or primitivity, in the case
that $r$ is a Mersenne exponent): just perform $r$ modular squarings,
starting from $x$, and check if the result is $x$. Since the cost of
each squaring is $O(r)$, the cost of the irreducibility test is $O(r^2)$.

There are more sophisticated algorithms for testing irreducibility, 
based on modular composition~\cite{KU08} and
fast matrix multiplication~\cite{rpb045}. However, these algorithms
are actually slower than the classical algorithm when applied to 
trinomials of degree less than about $10^7$.

When searching for irreducible trinomials of degree $r$, we can assume
that $s \le r/2$, since $x^r + x^s + 1$ is irreducible iff 
the reciprocal polynomial
$x^r + x^{r-s} + 1$ is irreducible. This simple observation saves
a factor of~$2$. In the following, we always assume that $s \le r/2$. 

\section{Degrees of Factors}

In order to predict the expected behaviour of our algorithm, we need to know
the expected distribution of degrees of irreducible factors. 
Our complexity estimates are based on the assumption that trinomials
of degree $r$ behave like the set of all polynomials of the same degree,
up to a constant factor:
\begin{hypothesis} \label{hyp1}
Over all trinomials $x^r+x^s+1$ of degree $r$ over $\GF(2)$,
the probability $\pi_d$ that a trinomial has no non-trivial factor of degree
$\le d$ is at most $c/d$, where $c$ is an absolute constant 
and $1 < d \le r/\ln r$.
\end{hypothesis}
This assumption is plausible and in agreement with experiments,
though not proven. It is not critical, because the correctness
of our algorithms does not depend on the assumption~-- only the predicted
running time depends on it. 
The upper bound $r/\ln r$ on $d$ is large enough for our application to
predicting the running time. An upper bound of $r$ on $d$ would probably
be incorrect, since it would imply at most $c$ irreducible trinomials
of degree $r$, but we expect this number to be unbounded.

Some evidence for the assumption, in the
case $r = r_{38}$, is presented in Table~2. The maximum
value of $d\pi_d$ is $2.08$, occurring at $d = 226\,887$. It would be 
interesting to try to explain the exact values of $d\pi_d$ for small $d$, 
but this would lead us too far afield.

\begin{table}[ht]
\begin{center}
\begin{tabular}{|c|c|}
\hline
 ~&\\[-12pt]
$d$ 	& $d\pi_d$	\\
\hline
1	& 1.00		\\
2	& 1.33		\\
3	& 1.43		\\
4	& 1.52		\\
5	& 1.54		\\
6	& 1.60		\\
7	& 1.60		\\
8	& 1.67		\\
9	& 1.64		\\
10	& 1.65		\\
100	& 1.77		\\
1000	& 1.76		\\
10000	& 1.88		\\
226887	& 2.08		\\
\hline
\end{tabular}
\end{center}
\caption{Statistics for $r = r_{38}$}
\label{tab:table2}
\end{table}

\section{Sieving}

When testing a large integer $N$ for primality, it is sensible 
to check if it has
any small factors before applying a primality test such as the
AKS, %
ECPP, %
or (if we are willing to accept a small
probability of error) Rabin-Miller %
test. %
Similarly, when testing a high-degree polynomial for irreducibility, it is
wise to check if it has any small factors before applying the $O(r^2)$ test.

Since the irreducible polynomials of degree $d$ divide ${\mathbf P}_d(x)$,
we can check if a trinomial $T$ has a factor of degree $d$ (or some
divisor of $d$) by computing
\[\gcd(T, {\mathbf P}_d).\]
If $T = x^r + x^s + 1$ and $2^d < r$, we can reduce this to the
computation of a GCD of polynomials of degree less than $2^d$. Let
$d' = 2^d-1$, $r' = r \bmod d'$, $s' = s \bmod d'$. Then
${\mathbf P}_d = x (x^{d'} - 1)$,
\[T = x^{r'} + x^{s'} + 1 \bmod (x^{d'}-1),\]
so we only need to compute
\[\gcd(x^{r'} + x^{s'} + 1, x^{d'}-1).\]
We call this process ``sieving'' by analogy with the process of sieving out
small prime factors of integers, even though it is performed using GCD
computations.

If the trinomials
that have factors of degree less than $\log_2(r)$ are excluded by sieving,
then by Assumption~\ref{hyp1} we are left with $O(r/\log r)$ trinomials to 
test. The cost of sieving is negligible.
Thus the overall search has cost $O(r^3/\log r)$.

\section{The Importance of Certificates}

Primitive trinomials of degree $r < r_{32} = 756\,839$ 
are listed in Heringa \emph{et al.}~\cite{Heringa}.
Kumada \emph{et al.}~\cite{Kumada00} reported a search
for primitive trinomials of degree %
$r_{33} = 859\,433$ (they did not consider $r_{32}$).
They found one primitive trinomial; however they missed the trinomial
$x^{859433} + x^{170340} + 1$, because of a bug in their sieving routine.
We discovered the missing trinomial in June 2000 while testing our program
on the known cases.

This motivated us to produce {\em certificates} of reducibility for all the
trinomials that we tested (excluding, of course, the small number that
turned out to be irreducible).  A certificate of reducibility is, ideally,
a non-trivial factor.  If a trinomial $T$
is found by sieving to have a small factor, then it is easy to keep a record
of this factor.  If we do not know a factor, but the trinomial fails the
irreducibility test~(\ref{eq:irred-cond}), then we can record the residue
$R(x) = x^{2^r} - x \bmod T$. Because the residue can be large, we might
choose to record only part of it, e.g., $R(x) \bmod x^{32}$.

\section{The Classical Period}

The period 2000--2003 could be called the {\em classical} period.
We used efficient implementations of the classical algorithms outlined
above.  Since different trinomials could be tested on different computers,
it was easy to conduct a search in parallel, using as many processors as
were available.  For example, we often made use of PCs in an undergraduate
teaching laboratory during the vacation, when the students were away.

In this way, we found 
three primitive trinomials of degree $r_{32} = 756\,839$ (in June 2000),
two of degree $r_{37} = 3\,021\,377$ (August and December 2000), 
and one of degree $r_{38} = 6\,972\,593$ (in August 2002)\footnote{Primitive 
trinomials of degree 
$r_{34}$, $r_{35}$ and $r_{36}$ were ruled out by Swan's theorem,
as were $r_{39}$ and $r'_{40}$.}. 
The computation for degree $r_{38}$ was completed and double-checked 
by July 2003. 

For degree $r_{38} = 6\,972\,593$, there turned out to be only one primitive
trinomial $x^r + x^s + 1$ (assuming, as usual, that $s \le r/2$)\footnote{
The unique primitive trinomial of degree $6\,972\,593$ is
$x^{6972593} + x^{3037958} + 1$. It was named {\em Bibury} after the village
that the three authors of~\cite{rpb214} were visiting on the day that it 
was discovered.}.
How can we
be sure that we did not miss any?  For each non-primitive trinomial we had a
certificate, and these certificates were checked in an independent
computation.  In fact, we found a small number of discrepancies, possibly
due to memory parity errors in some of the older PCs that were used. 
This is a risk in any long computation~-- we should not assume that
computers are infallible.
The same phenomenon was observed by Nicely~\cite{Nicely} in his computation
of Brun's constant (which also uncovered the infamous ``Pentium bug'').

Since we had caught up with the GIMPS project, we thought (not for the last
time) that this game had finished, and published our results
in~\cite{rpb199,rpb214}. However, GIMPS soon overtook us by finding several
larger Mersenne primes with exponents $\pm 1 \bmod 8$:
$r'_{41} = 24\,036\,583, \ldots, r'_{44} = 32\,582\,657$.

The search for degree $r_{38} = 6\,972\,593$ 
had taken more than two years
(February 2001 to July 2003), so
it did not seem feasible to tackle the new Mersenne exponents
$r'_{41}, \ldots, r'_{44}$.

\section{The Modern Period}

We realised that, in order to extend the computation, we had to find more
efficient algorithms. The expensive part of the computation was testing
irreducibility using equation~(\ref{eq:irred-cond}).  If we could sieve
much further, we could avoid most of the irreducibility tests.  From 
Assumption~\ref{hyp1}, if we could sieve to degree $r/\ln r$, then we 
would expect only $O(\log r)$ irreducibility tests.

What we needed was an algorithm that would find the smallest factor of
a sparse polynomial (specifically, a trinomial) in a time that was fast 
{\em on average}.

There are many algorithms for factoring polynomials over finite fields,
see for example \cite{vzGaGe99}. %
The cost of most of them is dominated by GCD computations.  However, it is
possible to replace most GCD computations by modular multiplications,
using a process called {\em blocking} (introduced by
Pollard~\cite{Pollard75} in the context of integer factorization,
and by von zur Gathen and Shoup~\cite{Gathen92} for polynomial
factorization).  The idea is simple: instead of computing
$\gcd(T,P_1), \ldots, \gcd(T,P_k)$ in the hope of finding
a non-trivial GCD (and hence a factor of $T$), 
we compute $\gcd(T, P_1 P_2 \cdots P_k \bmod T)$,
and backtrack 
if necessary to split factors if they are not irreducible.
Since a GCD typically takes about $40$ times as long as a modular
multiplication for $r \approx r'_{41}$, blocking can give a large
speedup.

During a visit by the second author to the first author in February 2007,
we realised that a second level of blocking could be used to replace
most modular multiplications by squarings. Since a modular multiplication
might take 400 times as long as a squaring 
(for $r \approx r'_{41}$), this second level of blocking can provide
another large speedup.  The details are described in~\cite{rpb230}.
Here we merely note that $m$ multiplications and $m$ squarings 
can be replaced by one
multiplication and $m^2$ squarings.  The optimal value of $m$ is
$m_0 \approx \sqrt{M(r)/S(r)}$, where $M(r)$ is the cost of a modular
multiplication and $S(r)$ is the cost of a modular squaring, and the
resulting speedup is about $m_0/2$. 
If $M(r)/S(r) = 400$, then $m_0 \approx 20$
and the speedup over single-level blocking is roughly a factor of ten.

Using these ideas, combined with a fast implementation of polynomial
multiplication (for details, see~\cite{rpb232}) and a
subquadratic GCD algorithm, we were able to find ten primitive trinomials
of degrees $r'_{41}, \ldots, r'_{44}$ by January 2008.
Once again, we thought we were finished and published our
results~\cite{rpb233}, only to have GIMPS leap ahead again
by discovering $M'_{45}$ in August 2008, and $M'_{46}$ and $M'_{47}$ shortly
afterwards. The exponent $r'_{46}$ was
ruled out by Swan's theorem, but we had to set to work on degrees
$r'_{45} = 43\,112\,609$ and (later) the slightly smaller
$r'_{47} = 42\,643\,801$. 

The search for degree $r'_{45}$ ran from September 2008 to May 2009, with
assistance from Dan Bernstein and Tanja Lange who kindly allowed us to use
their computing resources in Eindhoven, and resulted in 
four primitive trinomials of record degree.

The search for degree $r'_{47}$ ran from June 2009 to August 2009, and found
five primitive trinomials. In this case we were lucky to have access to a
new computing cluster with $224$ processors at the Australian National
University, so the computation took less time than the earlier searches.
 
The results of our computations in the ``Modern Period'' are given in 
Table~3. There does not seem to be any predictable
pattern in the $s$ values. The number of primitive trinomials for a given
Mersenne exponent $r = \pm 1 \bmod 8$ appears to follow a Poisson
distribution with mean about $3.2$ (and hence it is unlikely to be bounded by
an absolute constant~-- see the discussion of Assumption~1 above).

\begin{table}[ht]
\begin{center}
\begin{small}
\begin{tabular}{|c|c|} \hline
$r$ & $s$	\\ \hline
$24\,036\,583$ & $8\,412\,642$, 		%
	     $8\,785\,528$\\	 		%
$25\,964\,951$ & $880\,890$, $4\,627\,670$, $4\,830\,131$, 
	     $6\,383\,880$\\ 			%
$30\,402\,457$ & $2\,162\,059$\\ 		%
$32\,582\,657$ & $5\,110\,722$,	   		%
	     $5\,552\,421$, $7\,545\,455$\\
$42\,643\,801$ & $55\,981$, $3\,706\,066$, $3\,896\,488$,\\
	       & $12\,899\,278$, $20\,150\,445$\\
$43\,112\,609$ & $3\,569\,337$, $4\,463\,337$, $17\,212\,521$, $21\,078\,848$\\
\hline
\end{tabular}
\end{small}
 ~\\[2pt]
\caption{Primitive trinomials $x^r + x^s + 1$
whose degree $r$ is a Mersenne exponent, for $s \le r/2$.}
\label{tab:primitives}
\end{center}
\label{tab:table3}
\end{table}

\section{The Modern Algorithm -- Some Details} 

To summarize the ``modern'' algorithm for finding primitive trinomials,
we improve on the classical algorithm by sieving much further to find a factor
of smallest degree, using a factoring algorithm based on fast
multiplication and two levels of blocking.
In the following paragraphs we give some details of the modern algorithm
and compare it with the classical algorithms.

Given a trinomial $T = x^r + x^s + 1$, we search for a factor of
smallest degree $d \le r/2$. (In fact, using Swan's theorem,
we can usually restrict the search to $d \leq r/3$, because we know that
the trinomial has an odd number of irreducible factors.)
If such a factor is found, we know that $T$ is reducible, so the program
outputs ``reducible'' and saves the factor for a certificate of reducibility.
The factor can be found by taking the GCD of $T$ and $x^{2^d} + x$;
if this GCD is non-trivial, then $T$ has at least one factor of
degree dividing $d$. If factors of degree smaller than $d$ have already been
ruled out, then the GCD only contains factors of degree $d$ (possibly a
product of several such factors).
This is known as \emph{distinct degree factorization} (DDF).

If the GCD has degree $\lambda d$ for $\lambda > 1$,
and one wants to split the product into $\lambda$ factors of degree $d$, then
an \emph{equal degree factorization} algorithm (EDF) is used.
If the EDF is necessary it is usually cheap, since the
total degree $\lambda d$ is usually small if $\lambda > 1$.

In this way we produce certificates of reducibility that consist just of a
non-trivial factor of smallest possible degree, and the lexicographically
least such factor if there are several\footnote{It is worth going to the 
trouble to find the lexicographically least factor, since this makes the
certificate unique and allows us to compare different versions of the 
program and locate bugs more easily than would otherwise be the case.}.
The certificates can be checked, for example with an independent program
using NTL~\cite{NTL}, much faster than the original computation
(typically in less than one hour for any of the degrees listed in 
Table~3).

For large $d$, when $2^d \gg r$, we do not compute $x^{2^d} + x$
itself, but its remainder, say $h$, modulo $T$. Indeed, 
$\gcd(T,x^{2^d}+x) = \gcd(T,h)$.
To compute $h$, we start from $x$, perform $d$ modular squarings, and add
$x$. In this way, we work with polynomials of degree less than $2r$.
Checking for factors of degree $d$ costs $d$ modular
squarings and one GCD.
Since we check potential degrees $d$ in ascending order,
$x^{2^d} \bmod T$ is computed from $x^{2^{d-1}} \bmod T$, which was obtained
at the previous step, with one extra
modular squaring. Thus, from Table~1, 
the cost per value of $d$ is 
$O(M(r) \log r)$.
However, this does not take into account the speedup due to blocking, 
discussed above.

The critical fact is that most trinomials have a small factor, so the search
runs fast on average. 

After searching unsuccessfully for factors of degree $d < 10^6$ say, 
we could switch to the classical irreducibility test~(\ref{eq:irred-cond}),
which is faster than factoring if the factor has degree greater than 
about $10^6$.
However, in that case our list of certificates would be incomplete. Since it
is rare to find a factor of degree greater than $10^6$, we let the program
run until it finds a factor or outputs ``irreducible''.  In the latter case,
of course, we can verify the result using the classical test. 
Of the certificates (smallest  irreducible factors) 
found during our searches, %
the largest is a factor 
$P(x) = x^{10199457} + x^{10199450} + \cdots + x^4 + x + 1$
of the trinomial
$x^{42643801} + x^{3562191} + 1$.
Note that, although the trinomial is sparse and has a compact
representation, the factor is dense and hence too large to present
here in full.

\section{Classical versus Modern}

For simplicity we use the $\Otilde$ notation which ignores $\log$
factors.
The ``classical'' algorithm takes an expected 
time $\Otilde(r^2)$ per trinomial, or $\Otilde(r^3)$ to cover all
trinomials of degree~$r$.

The ``modern'' algorithm takes expected time $\Otilde(r)$ per trinomial,
or $\Otilde(r^2)$ to cover all trinomials of degree~$r$.

In practice, the modern algorithm is faster by a factor of about $160$ for
$r = r_{38} = 6\,972\,593$, 
and by a factor of about $1000$ for $r = r'_{45} = 43\,112\,609$.

Thus, comparing the computation for $r = r'_{45}$ with that
for $r = r_{38}$: using the classical algorithm would
take about $240$ times longer (impractical),
but using the modern algorithm saves a factor of $1000$.

\section{How to Speed up the Search}

The key ideas are summarised here. Points (1)--(4) apply to both the
classical and modern algorithms; points (5)--(6) apply only to the modern
algorithm.

\begin{enumerate}
\item Since the computations for each trinomial can be performed
independently, it is easy to conduct a search in parallel, using as many
computers as are available.
\item Because the coefficients of polynomials over $\GF(2)$ are just $0$ or
$1$, there is a one-one correspondence between polynomials of degree $<d$
and binary numbers with $d$ bits.  Thus, on a $64$-bit computer we can
encode a polynomial of degree $d$ in
$\lceil (d+1)/64\rceil$ computer words.
If we take care writing the programs, we can operate on such polynomials
using full-word computer operations, thus doing $64$ operations in parallel.
\item Squaring of polynomials over $\GF(2)$ can be done in {\em linear time}
(linear in the degree of the polynomial), because the cross terms in the
square vanish:
\[\left(\sum_k a_k x^k\right)^2 = \sum_k a_k x^{2k}\,. \]
\item Reduction of a polynomial of degree $2(r-1)$  modulo a trinomial 
$T = x^r + x^s + 1$ of degree $r$ can also be done in linear time.
Simply use the identity 
$x^n = x^{n+s-r} + x^{n-r} \bmod T$
for $n = 2r-2, 2r-3, \ldots, r$ to replace the terms of degree $\ge r$
by lower-degree terms.
\item Most GCD computations involving polynomials 
can be replaced by multiplication of polynomials, using a 
technique known as ``blocking'' (described above).
\item Most multiplications of polynomials can be replaced by squarings,
using another level of blocking, as described in~\cite{rpb230}.
\end{enumerate}

\section{Conclusion}

The combination of these six ideas makes it feasible to find
primitive trinomials of very large degree. 
In fact, the current record degree
is the same as the largest known Mersenne exponent, 
$r = r'_{45} = 43\,112\,609$.
We are ready to find more primitive trinomials as soon as GIMPS finds
another Mersenne prime that is not ruled out by Swan's Theorem.
Our task is easier than that of GIMPS, because finding a primitive trinomial
of degree $r$, and verifying that a single value of $r$ 
is a Mersenne exponent, both cost about the same: $\Otilde(r^2)$.

The trinomial hunt has resulted in improved software for
operations on polynomials over $\GF(2)$, and has shown that the best
algorithms in theory are not always the best in practice. 
It has also provided a large database of factors of trinomials over
$\GF(2)$, leading to several interesting conjectures which are a topic
for future research.

\subsection*{Acknowledgements}

We thank Allan Steel for verifying many of our primitive trinomials using 
Magma~\cite{Magma}, and 
Philippe Falandry, 
Shuhong Gao,
Robert Hedges, 
Samuli Larvala, 
Brendan McKay, 
{\'Eric} Schost, 
Julian Seward, 
Victor Shoup,
Andrew Tridgell and 
George Woltman for their advice and assistance in various ways.
Nate Begeman, 
Dan Bernstein,
Nicolas Daminelli,
Tanja Lange,
Ernst Mayer, 
Barry Mead, 
Mark Rodenkirch,
Juan Luis Varona, and 
Mike Yoder contributed machine cycles to the search.
Finally, we thank the University of Oxford, the Australian National University,
and INRIA for use of their computing facilities, and the Australian Research
Council for its support.


\begin{thebibliography}{99}

\bibitem{Magma}
W. Bosma and J. Cannon,
{\em Handbook of Magma Functions},
School of Mathematics and Statistics, University of Sydney, 1995.
\url{http://magma.maths.usyd.edu.au/}


\bibitem{rpb232}
R. P. Brent, P. Gaudry, E. Thom\'e and P. Zimmermann,
Faster multiplication in ${\mathrm GF}(2)[x]$,  
{\em Proc.\ ANTS VIII 2008}, 
{\em Lecture Notes in Computer Science} {\bf{5011}}, 153--166.
\url{http://wwwmaths.anu.edu.au/~brent/pub/pub232.html}

\bibitem{rpb045}
R. P. Brent and H. T. Kung, 
Fast algorithms for manipulating formal power series,
{\em J.\ ACM} {\bf{25}} (1978), 581--595.
\url{http://wwwmaths.anu.edu.au/~brent/pub/pub045.html}

\bibitem{rpb199}
R. P. Brent, S. Larvala and P.~Zimmermann,
A fast algorithm for testing reducibility of trinomials mod~2
and some new primitive trinomials of degree 3021377,
{\em Math.\ Comp.} {\bf{72}} (2003), 1443--1452.
\url{http://wwwmaths.anu.edu.au/~brent/pub/pub199.html}

\bibitem{rpb214}
R. P. Brent, S. Larvala and P.~Zimmermann, 
A primitive trinomial of degree 6972593,  
{\em Math.\ Comp.\ }{\bf{74}} (2005), 1001--1002,
\url{http://wwwmaths.anu.edu.au/~brent/pub/pub214.html}

\bibitem{rpb230}
R. P. Brent and P. Zimmermann,
A multi-level blocking distinct-degree factorization algorithm,
{\em Finite Fields and Applications:
Contemporary Mathematics} {\bf 461} (2008), 47--58.
arXiv:0710.4410v1,
\url{http://wwwmaths.anu.edu.au/~brent/pub/pub230.html}

\bibitem{rpb233}
R. P. Brent and P. Zimmermann,
Ten new primitive binary trinomials,
{\em Math.\ Comp.} {\bf 78} (2009), 1197--1199.
\url{http://wwwmaths.anu.edu.au/~brent/pub/pub233.html}

\bibitem{vzGaGe99}
J. von zur Gathen and J. Gerhard,
{\em Modern Computer Algebra},
Cambridge Univ.\ Press, 1999.

\bibitem{Gathen92}
J. von~zur Gathen and V. Shoup,
Computing Frobenius maps and factoring polynomials,
{\em Computational Complexity }{\bf{2}} (1992), 187--224.

\bibitem{Heringa}
J.~R.~Heringa, H.~W.~J.~Bl\"ote and A.~Compagner,
{New primitive trinomials of Mersenne-exponent degrees
for random-number generation,}
{\em International J.\ of Modern Physics C} {\bf\ 3} (1992),
561--564.

\bibitem{KU08}
K. Kedlaya and C. Umans,
Fast modular composition in any characteristic,
{\em Proc.\ FOCS 2008}, 146--155.

\bibitem{Kumada00}
T. Kumada, H. Leeb, Y. Kurita and M. Matsumoto,
New primitive $t$-nomials $(t = 3$,~$5)$ over $\GF(2)$
whose degree is a Mersenne exponent,
{\em Math.\ Comp.\ }{\bf{69}} (2000), 811--814.
Corrigenda: {\em ibid} {\bf{71}} (2002), 1337--1338.

\bibitem{Nicely}
T.~Nicely,
A new error analysis for Brun's constant,
{\em Virginia Journal of Science} {\bf{52}} (2001), 45--55.

\bibitem{Pellet78}
A.-E.~Pellet,
Sur la d\'ecomposition d'une fonction enti\`ere en facteurs
irr\'eductibles suivant un module premier $p$,
{\em Comptes Rendus de l'Acad\'emie des Sciences Paris}
{\bf{86}} (1878), 1071--1072.

\bibitem{Pollard75}
J. M. Pollard. 
A Monte Carlo method for factorization, 
{\em BIT} {\bf 15} (1975), 331--334,

\bibitem{Schonhage77}
A. Sch{\"o}nhage, 
Schnelle Multiplikation von Polynomen {\"u}ber K{\"o}rpern der 
Charakteristik $2$, 
{\em Acta Informatica} {\bf 7} (1977), 395--398.

\bibitem{NTL}
V.~Shoup,    
{NTL: A library for doing number theory}.
\url{http:www.shoup.net/ntl/}

\bibitem{Stickelberger97}
L.~Stickelberger,
\"Uber eine neue Eigenschaft der Diskriminanten algebraischer Zahlk\"orper,
{\em Verhandlungen des ersten Internationalen Mathematiker-Kongresses},
Z\"urich, 1897, 182--193.

\bibitem{Swan62}
R.~G.~Swan,
{Factorization of polynomials over finite fields,}
{\em Pacific J.~Math. }{\bf{12}} (1962), 1099--1106. 

\bibitem{Cunningham}
S.~Wagstaff, Jr.,
The Cunningham Project.
\url{http://homes.cerias.purdue.edu/~ssw/cun/}

\bibitem{GIMPS}
G.~Woltman {\em et al.},
GIMPS, The Great Internet Mersenne Prime Search.
\url{http://www.mersenne.org/}

\end{thebibliography}
\end{document}